\documentclass[12pt]{amsart}
\usepackage{amssymb}

\newcommand{\<}{\langle}
\renewcommand{\>}{\rangle}

\renewcommand{\a}{\alpha}
\renewcommand{\b}{\beta}
\renewcommand{\d}{\delta}
\newcommand{\D}{\Delta}
\newcommand{\e}{\varepsilon}
\newcommand{\g}{\gamma}

\renewcommand{\l}{\lambda}

\newcommand{\s}{\sigma}

\renewcommand{\th}{\theta}

\newcommand{\z}{\zeta}

\renewcommand{\span}{\mathop{\mathrm{span}}\nolimits}
\newcommand{\Ad}{\mathop{\mathrm{Ad}}\nolimits}
\newcommand{\dist}{\mathop{\mathrm{dist}}\nolimits}

\newtheorem{theorem}{Theorem}[section]
\newtheorem{lemma}{Lemma}[section]
\newtheorem{corollary}{Corollary}[section]
\newtheorem{proposition}{Proposition}[section]
\newtheorem{notation}{Notation}[section]

\begin{document}

\title[Matrix Algebras Converge]{Matrix Algebras Converge to the 
Sphere for Quantum 
Gromov--Hausdorff Distance}
\author{Marc A. Rieffel}
\address{Department of Mathematics \\
University of California \\ Berkeley, CA 94720-3840}
\email{rieffel@math.berkeley.edu}
\date{December 11, 2002}
\thanks{The research reported here was
supported in part by National Science Foundation grant DMS99-70509.}

\begin{abstract}
On looking at the literature associated with string theory one finds
statements that a sequence of matrix algebras converges to the $2$-sphere
(or to other spaces).  There is often careful bookkeeping with lengths,
which suggests that one is dealing with ``quantum metric spaces''.  We
show how to make these ideas precise by means of Berezin quantization
using coherent states.  We work in the general setting of integral
coadjoint orbits for compact Lie groups.
\end{abstract}

\maketitle
\allowdisplaybreaks

On perusing the theoretical physics literature which deals with string
theory and related parts of quantum field theory, one finds in many
scattered places assertions that the complex matrix algebras, $M_n$,
converge to the two-sphere, $S^2$, (or to related spaces) as $n$ goes to
infinity.  Here $S^2$ is viewed as synonymous with the algebra $C(S^2)$ of
continuous complex-valued functions on $S^2$ (of which $S^2$ is the
maximal-ideal space).  Approximating the sphere by matrix algebras is
attractive for the following reason.  In trying to carry out quantum field
theory on $S^2$ it is natural to try to proceed by approximating $S^2$ by
finite spaces.  But ``lattice'' approximations coming from choosing a
finite set of points in $S^2$ break the very important symmetry of the
action of $SU(2)$ on $S^2$ (via $SO(3)$).  But $SU(2)$ acts naturally on
the matrix algebras, in a way coherent with its action on $S^2$, as we
will recall below.  So it is natural to use them to approximate $C(S^2)$.  In
this setting the matrix algebras are often referred to as ``fuzzy
spheres''.  (See \cite{M2}, \cite{M}, \cite{GK2}, \cite{HL}, \cite{IKT}
and references therein.)

When using the approximation of $S^2$ by matrix algebras, the precise
sense of convergence is usually not explicitly specified in the
literature.  Much of the literature is at a largely algebraic level, with
indications that the notion of convergence which is intended involves how
structure constants and important formulas change as $n$ grows.  See, for
example, \cite{M3}, \cite{GP}, \cite{BMO} and section~IIC of \cite{Ta2}.
There is some discussion of approximation by matrix algebras at a more
analytical level within the mathematical 
physics literature concerned with quantization of symplectic
manifolds.  See the references in \cite{L}, \cite{Ln2}.  Much of this
discussion goes in the direction of showing that the matrix algebras can
be combined with $C(S^2)$ to form a continuous field of $C^*$-algebras
with $C(S^2)$ as limit point as $n$ grows.  To me the most satisfying
version of this idea has been given by Landsman \cite{Ln2} (or \cite{L}).
He shows, in the more general setting of coadjoint orbits of compact Lie
groups, and by means of Berezin quantization, that suitable matrix
algebras form a strict quantization, as defined in \cite{L}, \cite{R8},
\cite{R9}, \cite{R10}. This means that not only does one have a continuous
field of $C^*$-algebras, but also that commutators in the $M_n$'s converge to
the Poisson bracket on $S^2$ in a precise analytical sense (so that one
has a good ``semi-classical limit''). A more complicated proof, in the
more general setting of compact K\"ahler manifolds, was given earlier
in \cite{BMS}.

But if one goes back to the string-theory literature, one sees that there
is much more in play than just the continuous-field aspect.  Almost always
there are various lengths involved, and the writers are often careful in
their bookkeeping with these lengths as $n$ grows.  This suggested to me
that one is dealing here with metric spaces in some quantum sense, and
with the convergence of quantum metric spaces.  Now the only notion of
convergence of classical compact metric spaces with which I am familiar
is that of Gromov--Hausdorff convergence.  With this in mind, I gave in
\cite{R6} a definition of what one might mean by a compact quantum metric
space, and what can be meant by quantum Gromov--Hausdorff convergence of
these compact quantum metric spaces.  (``Compact'' here means that we
restrict attention to initial algebras.)  I showed that these notions have
basic properties which closely parallel the classical theory.

The purpose of the present paper is to show that when matrix algebras
are equipped in a natural way with a ``metric'', then they converge in
quantum Gromov--Hausdorff distance to $S^2$ with its round metric (for a
given radius).  In fact, we show that the corresponding fact is equally
well true for any integral coadjoint orbit (corresponding to an
irreducible representation) of a compact Lie group, once suitable
definitions are given.  (Examples in the field-theory and
string-theory literature of the approximation of coadjoint orbits other
than the two-sphere by matrix algebras are given in \cite{GS}, \cite{TrV},
\cite{ABI}, \cite{BDL}, \cite{DJ} and the references therein.)  
As in the work of Landsman, our
principal tool is Berezin quantization.

But when one goes back yet again to the 
field and string theory literature, one is
reminded that there is still much more structure involved than just a
metric structure and possible Gromov--Hausdorff convergence.  The field and
string-theorists need the whole apparatus of non-commutative differential
geometry --- ``vector bundles'', connections, ``Riemannian metrics'' --- so
that they can define action functionals such as Yang--Mills functionals.
(See \cite{GRS}, \cite{IKT} and references therein.)  Furthermore they
also want supersymmetry.  Thus what seems to be needed is a
``Gromov--Hausdorff'' convergence which encompasses this whole apparatus.
But so far as I am aware, such a theory does not yet exist even for
ordinary spaces.  But hints of what such a theory might look like can be
found in the literature concerned with the collapsing of Riemannian
manifolds.  The limit spaces need not be manifolds.  But certain parts of
the metric and differential geometric apparatus persist \cite{Fu},
\cite{G2}, \cite{KSo}, \cite{Sa}, \cite{Sh}, \cite{Lt1}, \cite{Lt2},
\cite{Sbl} \cite{Gro}.  It is an interesting challenge to give a full
characterization of the parts of the apparatus which do persist, and to
give an effective definition of the convergence of this apparatus; and
then to find the appropriate quantum generalizations.

An interesting aspect of our theory is that if ${\mathcal O}$ is a
coadjoint orbit for some compact Lie group, then there will be a suitable
sequence $\{n_j\}$ of dimensions such that the sequence of matrix algebras
$M_{n_j}$ converge to ${\mathcal O}$.  But this same sequence of matrix
algebras also converges to the sphere $S^2$, which need not be homeomorphic
to ${\mathcal O}$.  What is making the crucial difference is that the
metric structures which are placed on the matrix algebras are different in
the two cases.  This phenomenon may perhaps have some relation to the
ideas of ``change of topology'' which one finds in some places in the
string-theory literature \cite{MS}, \cite{Blc}, \cite{BBC}, 
\cite{Hyk}, \cite{CMS}. More specifically, it will follow from what 
we do that for any $\epsilon > 0$ there is a finite sequence of compact 
quantum metric spaces such that the first one is $\mathcal O$, the
last one is $S^2$, the intermediate ones are all full matrix
algebras, and the sequence is an $\epsilon$-chain in the sense that the 
quantum Gromov-Hausdorff distance between any two successive elements
of the chain is no greater than $\epsilon$.

The field and 
string theorists are interested in approximating many more spaces by
matrix algebras than just those which occur as coadjoint orbits.  For
example, they consider tori (\cite{Dj}, \cite{AM} and references therein),
surfaces of higher genus \cite{CDS}, \cite{BM}, and higher-dimensional
spheres \cite{HL}, \cite{Rmg}.  I am optimistic that the main theorem of
the present paper can in some form be extended (by different methods) to
more general homogeneous spaces, such as higher-dimensional spheres
\cite{Rmg}, since one still has available the action of a compact Lie
group which is so heavily used here.  I am also optimistic that tori will
not be difficult to treat, by using yet other methods of harmonic
analysis.  The substantial literature on quantization of K\"ahler
manifolds using Berezin quantization suggests that they too may well be
approximable in quantum Gromov--Hausdorff distance by matrix algebras with
suitable metrices.  (See \cite{Sch} and the references therein.)  This
should be an interesting project. My doctoral student Hanfeng Li has checked
that by using some of the facts in \cite{KrS2} one can carry out a certain
number of steps parallel to those which we carry out here for
coadjoint orbits. (One should also consider almost
K\"ahler manifolds \cite{KrS}.)  But it is far from clear to me how far
the theory might extend to arbitrary compact Riemannian manifolds.

As discussed in \cite{R5} and \cite{R6}, the appropriate way to specify a
``metric'' on a unital $C^*$-algebra, $A$, is by means of a seminorm
(generally unbounded) which plays the role of the Lipschitz seminorm on the
functions on an ordinary metric space.  In the present paper these
seminorms are defined in a quite simple way.  Let $G$ be a compact Lie
group, and let $\ell$ be a continuous length function on $G$ (for example,
coming from the Riemannian metric on $G$ corresponding to an Ad-invariant
inner product on the Lie algebra of $G$, especially if one wants the usual
round metric on the sphere).  Let $\a$ be an action of $G$ on 
$A$, and assume that the action is ``ergodic'' in the sense
that the only $\a$-invariant elements of $A$ are the scalar multiples of
the identity element of $A$.  Then we define the corresponding Lipschitz
seminorm, $L$, on $A$ by
$$
L(a) = \sup\{\|\a_x(a) - a\|/l(x): x \ne e\}, \leqno(0.1)
$$
where $e$ denotes the identity element of $G$.  (We may well have $L(a) =
+\infty$, but the set of $a$'s for which $L(a) < \infty$ is a dense
$*$-subalgebra.)  Some of the attractive properties of this definition
were discussed in \cite{R4}.  One pertinent example is the action of
$SO(3)$, and thus $SU(2)$, on $S^2$; but an equally pertinent example comes
from considering an irreducible unitary representation, $U$, of $G$ on a
finite-dimensional Hilbert space ${\mathcal H}$.  Let $A$ be $B({\mathcal
H})$, the algebra of operators on ${\mathcal H}$ (a full matrix algebra),
and let $\a$ be the action of $G$ on $A$ by conjugation by $U$.  Then
$(B({\mathcal H}),L)$ is a fine example of a compact quantum metric space
(and $L$ depends strongly on which group, with what representation, acts on
${\mathcal H}$).  For $G = SU(2)$ we will show that $(B({\mathcal H}),L)$
converges to $(C(S^2),L)$ for quantum Gromov--Hausdorff distance as the
dimension of ${\mathcal H}$ increases, and similarly for other coadjoint
orbits.

As mentioned above, our main analytical tool is Berezin quantization in
terms of coherent states.  What we need will be reviewed below.  But, with
$G = SU(2)$ and with $(U,H)$ as above, and for a choice of highest
weight-vector, $\xi$, in ${\mathcal H}$, Berezin defines for each $T \in
B({\mathcal H})$ its symbol, $\s_T$, which is a continuous function on the
coadjoint orbit, ${\mathcal O}$, of the weight for $\xi$.  When
$B({\mathcal H})$ is equipped with its Hilbert--Schmidt norm, and
$C({\mathcal O})$ with its $L^2({\mathcal O})$-norm, $\s$ has an adjoint,
${\breve \s}$, from $C({\mathcal O})$ to $B({\mathcal H})$.  The maps $\s$
and ${\breve \s}$ provide our principal tools for estimating the quantum
Gromov--Hausdorff distance from $B({\mathcal H})$ to $C({\mathcal O})$.
The composition $\s \circ {\breve \s}$ is called the Berezin transform, and
has received considerable study.  Nevertheless, the main estimate which we
need for the Berezin transform seems to be new, and of independent
interest.  We give a slightly imprecise statement of it here.  Label
${\mathcal H}$ by its dimension, thus ${\mathcal H}_n$.  Then $\s$ and
${\breve \s}$ depend on $n$, and we write $(\s \circ {\breve \s})_n$ for the
corresponding transform.

\setcounter{section}{3}
\setcounter{theorem}{3}
\begin{theorem}[imprecise]
\label{th3.4a}
There is a sequence $\{\d_n\}$ of numbers converging to $0$ such that
\[
\|f - (\s \circ {\breve \s})_n(f)\|_{\infty} \le \d_n L(f)
\]
for every $f \in C({\mathcal O})$ and every $n$.
\end{theorem}

I have not seen a relation such as this one between the Berezin transform
and Lipschitz norms discussed in the literature.

We also need information about the composition in the opposite order,
$({\breve \s} \circ \s)_n$, which carries $B({\mathcal H}_n)$ into itself.
I have not seen this mapping discussed in the literature at all.  We will
prove:

\setcounter{section}{6}
\setcounter{theorem}{0}
\begin{theorem}[imprecise]
\label{th6.1a}
Let $\g_n$ be the smallest constant such that
\[
\|T - ({\breve \s} \circ \s)(T)\| \le \g_n L(T)
\]
for all $T \in B({\mathcal H}_n)$.  Then the sequence $\{\g_n\}$ converges
to $0$.
\end{theorem}

The proof of this theorem involves other interesting facts about Berezin
symbols.

This paper is organized as follows.  In Section~\ref{sec1} we introduce
much of our notation and many of the structures which we need.  In
Section~\ref{sec2} we carry the development as far as we can for general
compact groups (including finite ones).  We turn to compact Lie groups in
Section~\ref{sec3}, where we state our main theorem (Theorem~\ref{th3.2})
concerning convergence of matrix algebras to coadjoint orbits.  We also
develop there the facts about Berezin covariant symbols which we need to
prove Theorem~\ref{th3.4a} (stated above).  Section~\ref{sec4} develops
further facts about covariant symbols, which are then used in
Section~\ref{sec5} to obtain the facts about Berezin contravariant
symbols which we need.  Then in Section~\ref{sec6} we prove
Theorem~\ref{th6.1a} (stated above), and use it to conclude the proof of
our main theorem.

Let me mention that David Kerr
has recently developed a matricial version of quantum Gromov-Hausdorff
distance \cite{Ker}, and he indicates how it applies to the topic
of the present paper.
Also, my former doctoral student Hanfeng Li 
has worked out \cite{Li} that the
main theorem (Theorem~\ref{th3.2}) can also be successfully approached 
by using the
results on continuous fields of quantum metric spaces developed in
\cite{R6} and the results of Landsman \cite{Ln2}, \cite{L} showing that
Berezin quantization gives a strict quantization. Here we will not
assume that the continuous-field structure is known (though
Theorem 4.2 and Proposition 4.3 provide the main steps in establishing
it). Thus our treatment here is more self-contained. 

A substantial part of the research reported here was carried out while I
visited the Institut de Math\'ematique de Luminy, Marseille, for three
months.  I would like to thank Gennady Kasparov, Etienne Blanchard, Antony
Wasserman, and Patrick Delorme very much for their warm hospitality and
their mathematical stimulation during my very enjoyable visit.


\setcounter{section}{0}
\section{The quantum metric spaces}
\label{sec1}

Let $G$ be a compact group (perhaps even finite, at first).  Let $U$ be an
irreducible unitary representation of $G$ on a Hilbert space ${\mathcal
H}$.  Let $B = B({\mathcal H})$ denote the $C^*$-algebra of all linear
operators on ${\mathcal H}$ (a ``full matrix algebra'', with its operator
norm).  There is a natural action, $\a$, of $G$ on $B$ by conjugation by
$U$.  That is, $\a_x(T) = U_xTU_x^*$ for $x \in G$ and $T \in B$.  Because
$U$ is irreducible, the action $\a$ is ``ergodic'', in the sense that the
only $\a$-invariant elements of $B$ are the scalar multiples of the
identity operator.  Fix a continuous length function, $\ell$, on $G$ (so
$G$ must be metrizable).  Thus $\ell$ is non-negative, $\ell(x) = 0$ iff
$x = e$ (the identity element of $G$), $\ell(x^{-1}) = \ell(x)$, and
$\ell(xy) \le \ell(x) + \ell(y)$.  For the reasons given in
Section~\ref{sec2} we also require a condition 
discussed in example 6.5 of \cite{R6},
namely that $\ell(xyx^{-1}) = \ell(x)$ for all $x$, $y \in G$.  Then in
terms of $\a$ and $\ell$ we can define a seminorm, $L_B$, on $B$ by the
formula $(0.1)$ given in the introduction.  Then $(B,L_B)$, or more
precisely its self-adjoint part, is an example of a compact quantum metric
space, as defined in \cite{R6}.  (For a quite different way of defining a
Lip-norm on $B$, for which the approximation of $A$ by $B$ in the quantum
Gromov--Hausdorff metric is almost a tautology, see equation $3.19$ of
\cite{ZS}.)

Let $P$ be a rank-one projection in $B({\mathcal H})$ (traditionally
specified by giving a non-zero vector in its range).  For any $T \in B$ we
define its Berezin covariant symbol \cite{Brz}, \cite{Prl}, $\s_T$, with
respect to $P$, by
\[
\s_T(x) = \tau(T\a_x(P)),
\]
where $\tau$ denotes the usual (un-normalized) trace on $B$.  (When the
$\a_x(P)$'s are viewed as giving states on $B$ via $\tau$ as above, they
form a family of ``coherent states'' \cite{Prl}.)  Let $H$ denote the
stability subgroup of $P$ for $\a$.  Then it is evident that $\s_T$ can be
viewed as a (continuous) function on $G/H$.  We let $\l$ denote the action
of $G$ on $G/H$, and so on $A = C(G/H)$, by left-translation.  If we 
note that
$\tau$ is $\a$-invariant, then it is easily seen that $\s$ is a unital,
positive, norm-nonincreasing, $\a$-$\l$-equivariant map from $B$ into $A$.

Of course, from $\l$ and $\ell$ we obtain a seminorm, $L_A$, on $A$ by
formula $(0.1)$.  It is just the restriction to $A$ of the seminorm on
$C(G)$ which we get from $\ell$ when we view $C(G/H)$ as a subalgebra of
$C(G)$, as we will often do when convenient.  We will often not restrict
$L_A$ to the Lipschitz functions, but rather permit $L_A(f) = +\infty$.  From
$L_A$ we obtain the usual quotient metric \cite{W2}
on $G/H$ coming from the metric
on $G$ for $\ell$.  One can check easily that $L_A$ in turn comes from
this quotient metric.  Thus $(A,L_A)$ is the compact quantum metric space
associated to this ordinary compact metric space.

It is reasonable to ask whether $\s$ might say something about the quantum
Gromov--Hausdorff distance (reviewed below) between $(A,L_A)$ and
$(B,L_B)$.  We stress that $H$, and so $A$, depends on the choice of the
projection $P$, and that even if two choices of $P$ have the same
stability group $H$, the symbol maps $\s$ may be different, and may give
quite different estimates of quantum Gromov--Hausdorff distance.
(Strictly speaking we should, according to \cite{R6}, be using the
self-adjoint parts of $A$ and $B$.  But by the comments just before
definition~$2.1$ of \cite{R6} we can, and will, be careless about this.)

According to the definition of quantum Gromov--Hausdorff distance given in
\cite{R6}, we must examine seminorms, $L$, on $A \oplus B$ whose quotient
seminorms on $A$ and $B$ are $L_A$ and $L_B$ respectively.  Furthermore,
$L$ must be a Lip-norm, as defined in \cite{R5}, meaning that the
null-space of $L$ is spanned by $(1_A,1_B)$, and that when $L$ is used to
define a metric, $\rho_L$, on the state space, $S(A \oplus B)$, of $A
\oplus B$, by
\[
\rho_L(\mu,\nu) = \sup\{|\mu(c)-\nu(c)|: c \in A \oplus B,\ L(c) \le 1\},
\]
then the topology on $S(A \oplus B)$ from $\rho_L$ coincides with the
weak-$*$ topology.  The state spaces $S(A)$ and $S(B)$ can be viewed in an
evident way as subsets of $S(A \oplus B)$.  We can then consider the usual
Hausdorff distance between them for $\rho_L$.  By definition \cite{R6},
the quantum Gromov--Hausdorff distance between $(A,L_A)$ and $(B,L_B)$ is
the infimum of these Hausdorff distances as $L$ varies.

We thus need a way to construct Lip-norms $L$ on $A \oplus B$.  As
discussed in section~5 of \cite{R6}, a convenient way to do this is to
look for seminorms $N$ on $A \oplus B$, called ``bridges'', and then set
\[
L(a,b) = L_A(a) \vee L_B(b) \vee N(a,b),
\]
where $\vee$ denotes ``maximum''.  The advantage of this is that $N$ can
be taken to be bounded.  For the present situation, we will take our
bridges $N$ to be of the very simple form
\[
N(f,T) = \g^{-1}\|f - \s_T\|_{\infty}
\]
for $f \in A$ and $T \in B$, where $\g$ is some positive constant.  This
constant must be taken large enough that the corresponding $L$ has $L_A$
and $L_B$ as quotient seminorms.  But for $B$ this is no problem.  Note
first:

\setcounter{section}{1}
\begin{proposition}
\label{prop1.1}
For any $T \in B$ we have
\[
L_A(\s_T) \le L_B(T).
\]
\end{proposition}

\begin{proof}
Since $\s$ is equivariant and does not increase norms, we
have
\begin{eqnarray*}
\|\lambda_x(\s_T) - \s_T\|_{\infty}/\ell(x) &=
&\|\s_{(\a_x(T)-T)}\|_{\infty}/\ell(x) \\
&\le &\|\a_x(T) - T\|/\ell(x) \le L_B(T)
\end{eqnarray*}
for every $x \in G$.
\end{proof}

\setcounter{corollary}{1}
\begin{corollary}
\label{cor1.2}
For $L$ defined from $N$ as above, the quotient of $L$ on $B$ is $L_B$,
regardless of the choice of $\g$.
\end{corollary}

\begin{proof}
It is clear from the definition of $L$ that the quotient of $L$ on $B$ is
no smaller than $L_B$.  But, given $T \in B$, we can simply take $f =
\s_T$.  From the above proposition we then have $L_A(f) \le L_B(T)$.  And
$N(f,T) = 0$, so that
\[
L(f,T) = L_A(f) \vee L_B(T) \vee \g^{-1}N(f,T) = L_B(T),
\]
as desired.
\end{proof}

Thus the difficult issue is how big $\g$ must be in order that the
quotient of $L$ on $A$ be $L_A$.

But assume that such a suitable $\g$ has been found.  For the
corresponding $L$ and $\rho_L$ we must estimate the Hausdorff distance
between $S(A)$ and $S(B)$.  Again, one half of this is quite simple.

\setcounter{proposition}{2}
\begin{proposition}
\label{prop1.3}
Let $\g$ be chosen such that the quotient of $L$ on $A$ is $L_A$.  Then
$S(A)$ is in the $\g$-neighborhood of $S(B)$ for $\rho_L$.
\end{proposition}

\begin{proof}
Let $\mu \in S(A) \subset S(A \oplus B)$.  We must find $\nu \in S(B)
\subset S(A \oplus B)$ such that $\rho_L(\mu,\nu) \le \g$.  A natural
guess for $\nu$ is $\nu = \mu \circ \s$.  Because $\s$ is positive and
unital, $\nu$ is indeed in $S(B)$.  We show that this choice of $\nu$
works.  Let $(f,T) \in A \oplus B$ with $L(f,T) \le 1$.  Then, in
particular, $\|f - \s_T\|_{\infty} \le \g$.  Thus
\begin{eqnarray*}
|\mu(f,T) - \nu(f,T)| &= &|\mu(f) - \nu(T)| = |\mu(f) - \mu(\s_T)| \\
&= &|\mu(f-\s_T)| \le \|f-\s_T\|_{\infty} \le \g.
\end{eqnarray*}
Since this is true for all such $(f,T)$, we have $\rho_L(\mu,\nu) \le \g$,
as desired.
\end{proof}

It is thus clear that to get the best estimates, we want $\g$ to be as
small as possible consistent with the quotient of $L$ on $A$ being $L_A$.

But we will still have the second difficult issue of showing that $S(B)$
is in a suitably small neighborhood of $S(A)$, so that their Hausdorff
distance is small.  It will take all of the discussion in
Sections~\ref{sec4} and \ref{sec5} for us to handle this issue for compact
semisimple Lie groups.  But we address first, in the next section, 
the first of these two difficult issues.

We remark that it might be interesting to see if the Stratonovich--Weyl
symbols of \cite{FGV} could be used in place of the Berezin symbol in the
definition of the bridge $N$.  However the Stratonovich symbol does not in
general carry positive operators to positive functions, and so the proof
of Proposition \ref{prop1.3} would not carry over directly.


\section{Choosing the bridge constant $\g$}
\label{sec2}

In this section we give one method for finding a $\g$ which will work for
the bridge $N$ of the previous section.  This method will be adequate for
our later purposes.  We continue with the same  notation as in the
previous section.

We choose for $G$ the Haar measure which gives $G$ total measure $1$; and
on $G/H$ we choose the image of this Haar measure, which is a
$G$-invariant measure giving $G/H$ total mass $1$.  This is conveniently
done by viewing $C(G/H)$ as a subalgebra of $C(G)$.  We will often not
distinguish between a point in $G$ and is image in $G/H$, with the context
making clear what is intended.  We will denote integration on $G$ or $G/H$
by $dx$, $dy$, etc.

As mentioned in the introduction, we put on $A = C(G/H)$ the inner product
from $L^2(G/H)$, while on $B = B({\mathcal H})$ we put its
Hilbert--Schmidt inner product, using now the normalized trace
$d^{-1}\tau$, where $d$ denotes the dimension of ${\mathcal H}$.  Then the
mapping $\s$ from $B$ to $A$ has an adjoint operator, ${\breve \s}$, from
$A$ to $B$.  For any $T \in B$, a function $f \in A$ such that ${\breve
\s}_F = T$ is called a Berezin contravariant symbol \cite{Brz}, \cite{Prl}
for $T$.  The mapping ${\breve \s}$ is often viewed as a quantization, since
it takes functions to operators.

We need a familiar formula for ${\breve \s}$.  For the reader's convenience,
and to put matters into our notation, we give the derivation of this
formula here.  For any $f \in A$ and $T \in B$ we have
\begin{eqnarray*}
d^{-1}\tau({\breve \s}_fT^*) &= &\<{\breve \s}_f,T\> = \<f,\s_T\> 
= \int f(x)(\s_T(x))^-dx \\
&= &\int f(x)\tau(\a_x(P)T^*)dx 
= \tau(\int f(x)a_x(P)dx\;T^*).
\end{eqnarray*}
Since this is true for all $T$, we obtain
\[
{\breve \s}_f = d \int f(x)\a_x(P)dx.
\]
It is well-known and easily seen that the ``second orthogonality
relation'' for irreducible representations of groups can be written as
\[
I = d \int \a_x(P)dx,
\]
where $I$ is the identity operator.  From this we see that ${\breve \s}$ is
unital.  (It was in order to make ${\breve \s}$ unital that we normalized
the trace.)  It is also evident that ${\breve \s}$ is positive,
norm-nonincreasing, and $\l$-$\a$-equivariant.  It is easy to see that
${\breve \s}$ corresponds to a ``pure-state quantization'' as discussed in
\cite{L}.

Our aim is to find a $\g$ such that the quotient on $A$
of our $L$, as defined in the
previous section, is $L_A$.  Now from the definition of $L$ it is
clear that the quotient is never less than $L_A$.  Thus what we 
would like to show
is that, given $f \in C(G/H)$, we can find $T \in B({\mathcal H})$ such
that $L(f,T) = L_A(f)$.  It is reasonable to try $T = {\breve \s}_f$.  Now
by the same argument as in the proof of Proposition~\ref{prop1.1}, using
the fact that ${\breve \s}$ is equivariant and norm-nonincreasing, we find
that $L_B({\breve \s}_f) \le L_A(f)$.  Thus it suffices to have
\[
L_A(f) \ge N(f,{\breve \s}_f) = \g^{-1}\|f - \s({\breve \s}_f)\|_{\infty}.
\]
That is, we seek $\g$ such that
\[
\|f - \s({\breve \s}_f)\|_{\infty} \le \g L_A(f)
\]
for all $f \in A$.  Now the mapping $f \mapsto \s({\breve \s}_f)$ has
received substantial study, and is often referred to as the Berezin
transform \cite{Sch}.  But I have not seen in the literature any
discussion of its relation to Lipschitz norms.  We need the following
frequently derived formula for the Berezin transform.
\begin{eqnarray*}
(\s({\breve \s}_f))(x) &= &\tau({\breve \s}_f\a_x(P)) 
= \tau\left(d \int
f(y)\a_y(P)dy\;\a_x(P)\right) \\
&= &d \int f(y)\tau(\a_y(P)\a_x(P))dy 
= d \int f(y)\tau(P\a_{y^{-1}x}(P))dy.
\end{eqnarray*}
The function $(y,x) \mapsto \tau(\a_y(P)\a_x(P))$ is important for the
theory, and is a ``transition probability'' on $G/H$ as discussed in
\cite{L}, \cite{Ln2}.

\begin{notation}
\label{note2.1}
{\em For any rank-one projection $P$ on ${\mathcal H}$ we define $h_P \in
C(G/H)$ by}
\[
h_P(x) = d\;\tau(P\a_x(P)).
\]
\end{notation}

It is easily checked that if $\xi$ is a vector of unit length in
the range of $P$, then $h_P(x) = d|\<U_x\xi, \xi\>|^2$.
We note that $h_P$ is non-negative, and that
\[
\int h(x)dx = \tau\left( P d\int \a_x(P)dx\right) = \tau(P)=1,
\]
so that $h_P(x)dx$ is a probability measure on $G/H$.  Furthermore,
\linebreak
$h_P(x^{-1}) = h_P(x)$ in the sense that $h_P(x^{-1}H) = h_P(xH)$ for $x
\in G$, because $\tau$ is $\a$-invariant.  Also, $h_P(e) = d$, so that
$h_P$ must be somewhat concentrated near $e$ $(= eH)$ in $G/H$.  The
formula we need then becomes
\setcounter{equation}{1}
\begin{equation}
\label{eq2.2}
(\s({\breve \s}_f))(x) = \int f(y)h_P(y^{-1}x)dy = \int f(xy^{-1})h_P(y)dy.
\end{equation}
We recognize this as an ordinary convolution (reflecting the
$\l$-equivariance of $\s \circ {\breve \s}$), if we view the functions as
defined in $G$ rather than $G/H$.

We must now bring the length function $\ell$ into the picture.  We have
used it to define the seminorm $L_A$ on $A = C(G/H)$, with its
corresponding metric, $\rho_A$, on $G/H$.  As mentioned earlier, $L_A$ is
then the Lipschitz seminorm for $\rho_A$.  This means that
\[
|f(x) - f(y)| \le L_A(f)\rho_A(x,y)
\]
for all $f \in C(G/H)$ and $x,y \in G/H$.  (We permit $L_A(f) = +\infty$.)  As
discussed in example $6.5$ of \cite{R6}, our extra requirement that
$\ell(zxz^{-1}) = \ell(x)$ implies that the action of $G$ on $A$ will
leave $L$ invariant.  Thus $G$ acts as isometries on $G/H$ for $\rho_A$,
that is, $\rho_A(zx,zy) = \rho_A(x,y)$ for all $z \in G$ and $x,y \in G/H$.

We are now ready to obtain the estimate which we need.  Let $f \in
C(G/H)$.  Because $h_P$ gives a probability measure, and $h_P(y^{-1}) =
h_P(y)$, we have
\begin{eqnarray*}
|f(x) &- &(\s({\breve \s}_f))(x)| = \left| \int (f(x) -
f(y))h_P(y^{-1}x)dy\right| \\
&\le &\int |f(x) - f(y)|h_P(y^{-1}x)dy 
\le L_A(f) \int
\rho_A(x,y)h_P(y^{-1}x)dy \\
&= &L_A(f) \int \rho_A (x,xy)h_P(y^{-1})dy 
= L_A(f) \int
\rho_A(e,y)h_P(y)dy.
\end{eqnarray*}
We have thus obtained:

\setcounter{theorem}{2}
\begin{theorem}
\label{th2.3}
With notation as above, we have
\[
\|f - \s({\breve \s}_f)\|_{\infty} \le \g L_A(f)
\]
for all $f \in A = C(G/H)$ if we set
\[
\g = \int_{G/H} \rho_A(e,y)h_P(y)dy.
\]
\end{theorem}

\setcounter{corollary}{3}
\begin{corollary}
\label{cor2.4}
For $\g$ chosen by the formula just above, the seminorm $L$ on $A \oplus
B$ defined by
\[
L(f,T) = L_A(f) \vee L_B(T) \vee \g^{-1} \|f - \s_T\|_{\infty}
\]
has $L_A$ as its quotient on $A$.
\end{corollary}

It follows from Proposition~\ref{prop1.3} that $S(A)$ is then in the
$\g$-neighborhood of $S(B)$ for $\rho_L$.  Consequently, in order to
obtain an upper bound on the quantum Gromov--Hausdorff distance between
$(A,L_A)$ and $(B,L_B)$, we must show also that $S(B)$ is in a suitably
small neighborhood of $S(A)$.  That is, for each $\nu \in S(B)$ we must
find $\mu \in S(A)$ close to $\nu$.  In view of our success
in Proposition 1.3, it is
reasonable to try setting $\mu = \nu \circ {\breve \s}$.  For this choice of
$\mu$ we can estimate $\rho_L(\mu,\nu)$ as follows.  Suppose that $(f,T)
\in A \oplus B$ and that $L(f,T) \le 1$, so that $\|f - \s_T\|_{\infty} \le
\g$.  Then, because ${\breve \s}$ is norm-non-increasing, we have
\begin{eqnarray*}
|\mu(f,T) - \nu(f,T)| &= &|\nu({\breve \s}_f - T)| \\
&\le &\|{\breve \s}_f - T\| \le \|{\breve \s}_f - 
{\breve \s}(\s_T)\| + \|{\breve
\s}(\s_T) - T\| \\
&\le &\|f - \s_T\|_{\infty} + \|{\breve \s}(\s_T) - T\| \le \g + \|{\breve
\s}(\s_T) - T\|.
\end{eqnarray*}
Thus any bound which we can obtain on $\|{\breve \s}(\s_T) - T\|$ for
$L_B(T) \le 1$ will give us a bound on the quantum Gromov--Hausdorff
distance.  I do not see an effective method for obtaining a good bound in
general, e.g., for finite groups.  But for compact semisimple Lie groups
we will see how to do this in Sections~\ref{sec4} and \ref{sec5}.

It should be stressed that $\g$ as above (as well as $\s$ and ${\breve \s}$
and $H$) depends on our choice of the rank-one projection $P$, as we will
glimpse further  in the next sections.  It is not clear to me how to make
an optimal choice of $P$ in general, e.g., for finite groups.


\section{Compact semisimple Lie groups}
\label{sec3}

There is a natural question which we did not address in the previous
section.  For $\s$ to be most useful, it is reasonable to expect that if
two operators have the same symbol then they should be the same operator, that
is, $\s$ should be faithful.  But if we consider the Hilbert--Schmidt
inner product on $B = B({\mathcal H})$ using $\tau$, it is clear that we
will have $\s_T = 0$ exactly if $T$ is orthogonal to the linear span of
$\{\a_x(P): x \in G\}$.  Thus we want this linear span to be all of $B$.
This will not happen in general, but I have found in the literature almost
no discussion of exactly when it does happen.  The only situation in
which it is well understood that the span is all of $B$ seems to be that in
which $G$ is a compact semisimple Lie group and $P$ corresponds to a
highest weight vector.  Most proofs for this case use the complex
structure which one then has on $G/H$ \cite{Prl}.  However, Simon
\cite{Smn} has given a relatively elementary proof.  (A more complicated
but elementary
proof appears in theorem $4.11$ of \cite{Wld}, while an earlier variant of
it appears as theorem~$1$ of \cite{FGV}.)  Simon also gives a few examples
showing that the span can be all of $B$ for some other weight vectors, but
not for all.  There is additional discussion of these matters in
\cite{Wld}, \cite{MAC}.  Since the case of a highest weight vector is
crucial for our purposes, we include a proof here, along the lines of
Simon's proof \cite{Smn}, but very slightly simpler.  This also gives us
an opportunity to introduce some of the notation which we will use
extensively later on. We remark that in what follows it is not 
essential that $G$ be semisimple --- we could permit it to be the 
product of a semisimple group with a torus. But this would provide
no new coadjoint orbits, and would somewhat complicate our notation.

Thus, let $G$ be a connected compact semisimple Lie group.  We denote its
Lie algebra by ${\mathfrak g}_0$, because we will usually work with the
complexification, ${\mathfrak g}$, of ${\mathfrak g}_0$.  We choose a
maximal torus in $G$, with corresponding Cartan subalgebra of ${\mathfrak
g}$, its set of roots $\D$, positive roots $\D^+$, and elements
$\{H_{\b},E_{\b},E_{-\b}: \b \in \D^+\}$ for ${\mathfrak g}$, where
$[E_{\b},E_{-\b}] = H_{\b}$, etc., much as in equations VIII.56 of
\cite{Sm1}.  Let $(U,{\mathcal H})$ be an irreducible unitary
representation of $G$.  We let $U$ also denote the corresponding
representation of ${\mathfrak g}$.  Then $(U,{\mathcal H})$ will have a
highest weight vector, $\xi$, of norm~$1$, unique up to a scalar multiple,
and characterized by the fact that $U_{E_{\b}}\xi = 0$ for all $\b \in
\D^+$.

\setcounter{theorem}{0}
\begin{theorem}
\label{th3.1}
Let $P$ be the rank-one projection which has the highest weight vector
$\xi$ in its range.  Then
\[
\span\{\a_x(P): x \in G\} = B({\mathcal H}).
\]
\end{theorem}

\begin{proof}
Let $S$ denote the linear span of the $\a_x(P)$'s.  It is clear that $S$
is an $\a$-invariant subspace of $B({\mathcal H})$ which contains $P$.  We
let $\a$ also denote the corresponding action of ${\mathfrak g}$ on
$B({\mathcal H})$, given by
\[
\a_X(T) = [U_X,T] = U_XT - TU_X
\]
for $X \in {\mathfrak g}$ and $T \in B({\mathcal H})$.  Then this action of
${\mathfrak g}$ carries $S$ into itself.  From equation VIII.5.6a of
\cite{Sm1} we know that we can choose the $E_{\b}$'s so  that
$(U_{E_{\b}})^* = U_{E_{-\b}}$.  For $\eta,\z \in {\mathcal H}$ let
$\<\eta,\z\>_{\mathcal K}$ denote the rank-one operator on ${\mathcal
H}$ defined by $\<\eta,\z\>_{\mathcal K} \th =
\<\th,\z\>\eta$.  Then for $\b \in \D^+$ we have
\begin{eqnarray*}
\a_{E_{\b}}(P) &= &\a_{E_{\b}}(\<\xi,\xi\>_{\mathcal K}) \\
&= &\<U_{E_{\b}}\xi,\xi\>_{\mathcal K} + \<\xi,U_{E_{-\b}}\xi\>_{\mathcal
K} = \<\xi,U_{E_{-\b}}\xi\>_{\mathcal K}.
\end{eqnarray*}
In the same way, if we apply a product $E_{\b_1} \dots E_{\b_k}$ to $P$, 
where $\b_j \in \D^+$ for each $j$, we
obtain $\<\xi,\eta\>_{\mathcal K}$, where
\[
\eta = U_{E_{-\b_k}} \dots U_{E_{-\b_1}}\xi.
\]
But it is known (proposition VII.2 of \cite{Srr}) that the various
$\eta$'s of this form span ${\mathcal H}$.  In this way we see that $S$
must contain all rank-one operators of the form $\<\xi,\eta\>_{\mathcal
K}$ for $\eta \in {\mathcal H}$.  But $S$ is invariant for the $G$-action
$\a$.  It follows that $S$ contains all rank-one operators of the form
$\<U_x\xi,\eta\>_{\mathcal K}$ for $x \in G$.  But the various $U_x\xi$'s
span ${\mathcal H}$.  From this it follows that $S$ contains all rank-one
operators, and so coincides with $B({\mathcal H})$.
\end{proof}

We will let $\Ad$ denote the adjoint action of $G$ on ${\mathfrak g}_0$ and
${\mathfrak g}$, and let $\Ad^*$ denote the coadjoint action of $G$ on the
dual vector spaces ${\mathfrak g}_0^*$ and ${\mathfrak g}^*$.  Let $P$ be
the rank-one projection for the highest weight vector $\xi$.  Define an
element $\omega$ of ${\mathfrak g}^*$ by
\[
\omega(X) = -i\tau(U_XP) = -i\<U_X\xi,\xi\>.
\]
Note that $\omega$ is real on ${\mathfrak g}_0$, so it can be viewed as
an element of ${\mathfrak g}_0^*$.  Now
\[
(\Ad_x^*\omega)(X) = \omega(\Ad_{x^{-1}}(X)) = -i\tau(U_{x^{-1}}U_XU_xP) =
-i\tau(U_X \a_x(P)).
\]
 From this we see that the stability subgroup for $\omega$ under $\Ad^*$
coincides with the stability subgroup, $H$, for $P$ under $\a$, as is
well-known \cite{Ln2}, \cite{L}.  The coadjoint orbit, ${\mathcal O}$, of
$\omega$ is naturally identified with $G/H$ via $x \mapsto 
\Ad_x^*(\omega)$.  It
is worth remarking that $H$ clearly contains the center of $G$, and so it
is not important here whether $G$ is simply connected.  But if, in fact,
$G$ is simply connected, then to every coadjoint orbit which is integral
in the usual algebraic sense there will correspond an irreducible
representation of $G$ with highest weight vector giving that orbit.  
Thus if we
need to apply our metric considerations to all integral coadjoint orbits,
we should take $G$ to be simply connected (or work with projective
representations). Of course, the coadjoint orbit for a highest weight 
vector for $SU(2)$ will be a two-sphere $S^2$.

We now fix the representation $(U,{\mathcal H})$, and thus the coadjoint
orbit ${\mathcal O}$.  We will follow a path established by Berezin
\cite{Brz}, and then followed by many others, including Simon \cite{Smn}.  
See \cite{Ln2}, \cite{L} for a nice account.  For any $n$ we can form
$(U^{\otimes n},{\mathcal H}^{\otimes n})$, the $n^{\rm th}$ inner tensor
power of $(U,{\mathcal H})$.  Let $(U^n,{\mathcal H}^n)$ denote the
subrepresentation generated by $\xi^n = \xi^{\otimes n}$.  It is
well-known (e.g., as a consequence of proposition VII.2 of \cite{Srr})
that $(U^n,{\mathcal H}^n)$ is irreducible, with $\xi^n$ as highest weight
vector.  The weight for $\xi^n$ is easily calculated to be $n\omega$.  Thus
the corresponding orbit is $n{\mathcal O}$, which we identify with
${\mathcal O}$ by dividing by $n$.  In particular, the stability subgroup
of $n\omega$ is clearly still $H$.  We let $B^n = B({\mathcal H}^n)$.  The
action of $G$ on $B^n$ by conjugation by $U^n$ will be denoted again
simply by $\a$.  We denote the corresponding Lip-norm (for $\ell$) on
$B^n$ by $L_n$, and we denote the Lip-norm on $A = C(G/H)$ still by $L_A$.

We are now prepared to state the main theorem of this paper.

\setcounter{theorem}{1}
\begin{theorem}
\label{th3.2}
With notation as above, the quantum metric spaces $(B^n,L_n)$ converge to
$(A,L_A)$ for quantum Gromov--Hausdorff distance as $n$ goes to $\infty$.
\end{theorem}

Our proof of this theorem will extend over the remaining sections of this
paper.  But we take the initial steps  here.  We let $P^n$ denote the
rank-one projection for $\xi^n$.  We denote the corresponding Berezin
symbol map from $B^n$ to $A = C(G/H) = C({\mathcal O})$ by $\s^n$.  Let
$d_n$ denote the dimension of ${\mathcal H}^n$, and let $h_n = h_{P^n}$ as
in Notation~\ref{note2.1}.  A crucial fact for our purposes is:

\setcounter{lemma}{2}
\begin{lemma}
\label{lem3.3}
The sequence of probability measures $h_n(x)dx$ converges in the weak-$*$
topology to the $\d$-function at $e$; that is, $\int_{G/H} f(x)h_n(x)dx$
converges to $f(e)$ for every $f \in C(G/H)$.
\end{lemma}

A proof of this fact, but in a slightly more complicated context, is given
in proposition~4b of \cite{Dff}.  That proof is not elementary, because it
depends on facts about the rate of growth of the dimensions of
representations.  We give here, for our context, a very elementary proof.
(We remark that this lemma corresponds to condition~2 of definition~3 of
\cite{Ln2}.)

\begin{proof}[Proof of Lemma \ref{lem3.3}]
For $x \in G/H$ we have
\begin{eqnarray*}
h_n(x) &= &d_n\tau(P^n\a_x(P^n)) = d_n|\<U_x^n\xi^n,\xi^n\>|^2 \\
&= &d_n|\<U_x\xi,\xi\>|^{2n} = d_n(\tau(P\a_x(P)))^n.
\end{eqnarray*}
Set $g(x) = \tau(P\a_x(P))$, so that $h_n = d_ng^n$, where $g^n$ now means
the $n^{\rm th}$ power of $g$.  Note that $g(e) = 1$.  Suppose that $x \in
G$ and $g(x) = 1$.  Then the Hilbert--Schmidt inner product tells us that
$\a_x(P)$ must agree with $P$.  But this says that $x \in H$.  In other
words, for $g$ viewed on $G/H$, the only point $x$ at which $g(x) = 1$ is
$x = e$.  It is clear that $0 \le g(x) \le 1$ for any $x$.  Notice that
$d_n = \|g^n\|_1^{-1}$, so that for the moment the $d_n$'s can be viewed
simply as normalizing constants to obtain probability measures.

Thus let $X$ be any compact space with distinguished point $e$ and full
measure $dx$.  Let $g \in C(X)$ with $0 \le g \le 1$, and with $g(x) = 1$
exactly if $x = e$.  We need to know that the probability measures
$g^n/\|g^n\|_1$ converge to the delta-function at $e$.  The simple
argument that this is true is given in the course of the proof of theorem
$8.2$ of \cite{R6}.
\end{proof}

As before, we assume that we have fixed a continuous length function,
$\ell$, on $G$, and that $\rho_A$ is the corresponding $G$-invariant
metric on $G/H$ as discussed in the previous section.  For each $n$ let
$\g_n$ be defined as in Theorem~\ref{th2.3}, that is,
\[
\g_n = \int_{G/H} \rho_A(e,y)h_n(y)dy.
\]
Since $\rho_A$ is continuous and $\rho_A(e,e) = 0$, and since $h_n$
converges to the $\d$-function at $e$, it is now clear that the sequence
$\{\g_n\}$ converges to $0$.  Putting all of the above together with
Theorem~\ref{th2.3}, we obtain the following fact about the Berezin
transform, which seems to be new and of independent interest.

\setcounter{theorem}{3}
\begin{theorem}
\label{th3.4}
Let $G$ be a connected
compact semisimple Lie group with length function $\ell$, 
and let $\omega$ be an element of
an integral coadjoint orbit for $G$.  Let $H$ be the stabilizer of $\omega$,
and let $L_A$ be the Lip-norm on $A = C(G/H)$ corresponding to $\ell$.
For each integer $n$ let $(U^n,{\mathcal H}^n)$ be the irreducible unitary
representation of $G$ with highest weight $n\omega$, and let $B^n =
B({\mathcal H}^n)$ with its action of $G$ by conjugation by $U^n$.  Let
$\s^n$ be the Berezin covariant symbol  map from $B^n$ to $A$ using the
projection on the highest weight vector for $n\omega$, and let ${\breve \s}^n$
be its adjoint.  Then there is a sequence of numbers, $\{\g_n\}$,
converging to $0$, such that
\[
\|f - \s^n({\breve \s}^n(f))\|_{\infty} \le \g_nL_A(f)
\]
for all $f \in A$ and all $n$.
\end{theorem}

By the comment right after Corollary~\ref{cor2.4} it follows that, in
terms of the metrics $\rho_L$ on $S(A \oplus B^n)$ defined there using $\s^n$
for the bridges for $L$, we can find for any $\e > 0$ an integer $N$ such that
$S(A)$ is in the $\e$-neighborhood of $S(B^n)$ for all $n \ge N$.  The
next sections are devoted to showing that for $N$ large enough it is also
true that $S(B^n)$ is in the $\e$-neighborhood of $S(A)$.


\section{Covariant symbols}
\label{sec4}

At the end of Section~\ref{sec2} we indicated the importance of studying
the mappings ${\breve \s}^n \circ \s^n$ on $B^n$.  This requires a more
careful study of $\s^n$ and of ${\breve \s}^n$.  This section is devoted to
obtaining the information which we need about $\s^n$.

We use the notation of the previous section.  Since ${\mathcal O} \subset
{\mathfrak g}_0^*$, we can view ${\mathfrak g}$ as the space of linear
complex-valued polynomials on ${\mathfrak g}_0^*$, which we can then
restrict to ${\mathcal O}$.  For $X \in {\mathfrak g}$ we define $\Phi_X$
on ${\mathcal O}$ by
\[
\Phi_X(\Ad_x^*\omega) = i(\Ad_x^*\omega)(X) = i\omega(\Ad_{x^{-1}}(X)) =
\<U_XU_x\xi,U_x\xi\>.
\]
We include the factor of $i$ because when the definition of $\s$ which we
gave in Section~\ref{sec1} is applied here, we see that
\[
\s(U_X)(x) = \tau(U_X\a_x(P)) = \<U_XU_x\xi,U_x\xi\> = \Phi_X(x).
\]
Here we are using the identification $x \mapsto \Ad_x^*\omega$ of $G/H$
with ${\mathcal O}$ to view $\Phi_X$ as a function on $G/H$.  Note also
that $\lambda_y\Phi_X = \Phi_{\Ad_y(X)}$ for $x \in G$ and $X \in {\mathfrak
g}$, so that $\Phi$ is equivariant.

We let ${\mathcal P}({\mathcal O})$ denote the (unital) algebra of
functions on ${\mathcal O}$, or $G/H$, generated by all the $\Phi_X$'s.
The algebra ${\mathcal P}({\mathcal O})$ is clearly carried into itself by
the action $\Ad^*$ of $G$ on ${\mathcal O}$.  We denote this action on
${\mathcal P}({\mathcal O})$ again by $\Ad$, or, when we view ${\mathcal
P}({\mathcal O})$ on $G/H$, by $\lambda$.

Let ${\mathcal T}$ denote the full tensor algebra over ${\mathfrak g}$,
whose homogeneous parts are the various tensor powers ${\mathfrak
g}^{\otimes m}$.  Again $G$ acts on ${\mathcal T}$ via $\Ad$ (the diagonal
action).  By the universal property of ${\mathcal T}$ we extend $\Phi$ to
an algebra homomorphism from ${\mathcal T}$ onto ${\mathcal P}({\mathcal
O})$ defined on elementary tensors by the product
\[
\Phi(X_1 \otimes \dots \otimes X_k) = \prod_j \Phi_{X_j}.
\]
It is easily seen that $\Phi$ is still $\Ad$-equivariant.  (All of this
is related, of course, to the universal enveloping algebra of
${\mathfrak g}$, but we don't need that structure here.)

We let $(U^n,{\mathcal H}^n,\xi^n)$ be as defined in the previous
section.  Then for $X \in {\mathfrak g}$ and for any $n$ we have
\begin{eqnarray*}
\s^n(U_X^n)(e) &= &\<U_X^n\xi^n,\xi^n\> \\
&= &\<U_X\xi \otimes \xi \dots \otimes \xi,\xi^n\> + \<\xi \otimes U_X\xi
\otimes \xi \dots,\xi^n\> + \dots \\
&= &n\s(U_X)(e).
\end{eqnarray*}
By the equivariance of $\s$ it follows that
\[
\s^n(U_X^n)(x) = n\s(U_X)(x)
\]
for every $x \in G$, that is, $\s^n(U_X^n) = n\s(U_X)$.  It is then
natural to define a linear map, $\Phi^n$, from ${\mathfrak g}$ to $B^n$ by
\[
\Phi^n(X) = n^{-1}U_X^n.
\]

\begin{notation}
\label{note4.1}
{\em We extend $\Phi^n$ to an (Ad-equivariant) 
algebra homomorphism from ${\mathcal T}$ to
$B^n$, defined on elementary tensors by
\[
\Phi^n(X_1 \otimes \dots \otimes X_k) = \Phi^n(X_1)\Phi^n(X_2)\dots
\Phi^n(X_k),
\]
where the order of the terms is now important.}
\end{notation}

The next theorem has been known, at least in part, for many special
cases.  It is stated as theorem~2A in \cite{Dff} without proof, but with
attribution to Gilmore \cite{Glm}.  But Gilmore does not give complete
details, and the details which he does give can be simplified
substantially.  Since this theorem is crucial for our purposes, we include
a proof here. We will let
${\mathcal T}^m$ denote the direct sum of the homogeneous subspaces
${\mathfrak g}^{\otimes k}$ for $k \le m$.

\setcounter{theorem}{1}
\begin{theorem}
\label{th4.2}
For any $Z \in {\mathcal T}$ the sequence $\{\s^n(\Phi_Z^n)\}$ converges
uniformly on $G/H$ to $\Phi_Z$.  For any integer $m$ this convergence is
also uniform in $Z$ as $Z$ ranges over any bounded subset of ${\mathcal
T}^m$.
\end{theorem}

\begin{proof}
It clearly suffices to prove the first part of the above theorem for
homogeneous $Z$'s, and, in fact, for elementary tensors.  So assume that
$Z = X_1 \otimes X_2 \otimes \dots \otimes X_k$.  Then
\[
\s^n(\Phi_Z^n)(e) = n^{-k}\<U_{X_1}^nU_{X_2}^n \dots
U_{X_k}^n\xi^n,\xi^n\>.
\]
Notice that
\[
U_{X_k}^n\xi^n = U_{X_k}\xi \otimes \dots \otimes \xi + \xi \otimes
U_{X_k}\xi \otimes \xi \dots \otimes \xi + \dots \quad \quad.
\]
 From this we see that when $U_{X_1}^n\dots U_{X_k}^n\xi^n$ is fully
expanded into elementary tensors we will have $n^k$ terms.  We now assume
that $n > k$.  Then the number of terms all of whose components are only
of the form either $\xi$ or $U_{X_j}\xi$ is $n(n-1)\dots (n-k+1) = n^k -
p(n)$, where $p$ is a polynomial in $n$ of degree $k-1$.  Each of these
terms will give a contribution toward $\s^n(\Phi_Z^n)(e)$ of the form
$\prod\<U_{X_j}\xi,\xi\> = \Phi_Z(e)$.  There will be $p(n)$ other terms.
Let
\[
K_Z = \max\{1,\|U_{X_j}\|: j = 1,\dots,k\}.
\]
Each of the $n^k$ terms has norm no bigger than $K_Z^k$.  Thus
\begin{eqnarray*}
|\Phi_Z(e) &- &\s^n(\Phi_Z^n)(e)| \\
&= &|\Phi_Z(e) -
n^{-k}((n^k-p(n))\Phi_Z(e) + (p(n) \mbox{ remaining terms}))| \\
&\le &n^{-k}2p(n)K_Z^k.
\end{eqnarray*}

Now $\Ad_x(Z)$ is also an elementary tensor, and $K_{\Ad_x(Z)} = K_Z$. Thus
the above discussion applies to $\Ad_x(Z)$ as well.  Since $\Phi_Z(x) =
(\lambda_x\Phi_Z)(e) = \Phi_{\Ad_x(Z)}(e)$, we obtain
\[
\|\Phi_Z - \s^n(\Phi_Z^n)\|_{\infty} \le n^{-k}2p(n)K_Z^k.
\]
This clearly goes to $0$ as $n \rightarrow \infty$ since $p$ is of degree
$k-1$.

By considering a basis for ${\mathcal T}^m$ (consisting of elementary
tensors if desired), the statement about the convergence being uniform on
bounded subsets of ${\mathcal T}^m$ follows easily.
\end{proof}

\setcounter{proposition}{2}
\begin{proposition}
\label{prop4.3}
$\Phi({\mathcal T})$ is uniformly dense in $C({\mathcal O})$.
\end{proposition}

\begin{proof}
Already $\Phi({\mathfrak g})$ alone separates the points of ${\mathcal
O}$, for if $x \in G$ is such that $\Phi_X(\Ad_x^*(\omega)) =
\Phi_X(\omega)$ for all $X \in {\mathfrak g}$, then $\Ad_x^*(\omega) =
\omega$, so that $x \in H$.  This suffices, by $\Ad$-equivariance.  Next,
$\Phi_X$ is pure imaginary for $X \in {\mathfrak g}_0$, so that
$\Phi({\mathcal T})$ is closed under complex conjugation.  By definition
$\Phi({\mathcal T})$ contains the constant functions.  We can thus apply
the Stone--Weierstrass theorem.
\end{proof}

We now prepare for some approximations which we will make in the next
sections.  We let ${\hat G}$ denote as usual the set of equivalence
classes of irreducible unitary representations of $G$.  Let $S$ be a
finite subset of ${\hat G}$, which we fix for the rest of this section.
We let $A_S$ denote the direct sum of the isotypic components of $A$ for
the representations in $S$ and the action $\lambda$.  We call it the
$S$-isotypic subspace of $A$.  In the same way we will speak of the
$S$-isotypic subspaces of other representations of $G$.

\setcounter{lemma}{3}
\begin{lemma}
\label{lem4.4}
The $S$-isotypic subspace $A_S$ is finite dimensional, 
and there is an integer $q$ $(= q_S)$ such
that $\Phi({\mathcal T}^q) \supseteq A_S$. (We fix this integer $q$, 
in terms of $S$, for
most of the remainder of this paper.)
\end{lemma}

\begin{proof}
From the Peter--Weyl theorem the various isotypic components of $C(G)$
for $\lambda$ are all finite dimensional. But $A = C(G/H)$ is
a $\lambda$-invariant subspace of $C(G)$,
and thus $A_S$ is
finite dimensional.

Because $\Phi({\mathcal T})$ is dense in $A$ by Proposition~\ref{prop4.3},
when we compose $\Phi$ with the usual projection onto the sum of the
isotypic components for $S$, this composition will carry ${\mathcal T}$
onto $A_S$.  From the equivariance of $\Phi$ it follows that
$\Phi({\mathcal T}) \supseteq A_S$.  
By choosing preimages for a basis for $A_S$
we find $q$ such that $\Phi({\mathcal T}^q) \supseteq A_S$.
\end{proof}

We remark that, conversely, given any $Z \in {\mathcal T}$, it is in some
${\mathcal T}^m$, which is finite dimensional and $\Ad$-invariant; and so
there is a finite $S \subseteq {\hat G}$ such that $\Phi_Z \in A_S$.

For the rest of this paper we fix an $\Ad$-invariant inner product on
${\mathfrak g}$.  It gives a corresponding inner product on each
${\mathfrak g}^{\otimes k}$, and so on ${\mathcal T}$, where we take the
${\mathfrak g}^{\otimes k}$'s to be orthogonal to each other for different
$k$'s.  Norms of elements of ${\mathcal T}$, and of operators from
${\mathcal T}$, will always be defined with respect to this inner product
(and, for operators, a specified norm on the target space).

\setcounter{lemma}{4}
\begin{lemma}
\label{lem4.5}
There is a constant, $K$, such that when we view each $\Phi^n$ as an
operator from ${\mathcal T}^q$ to $B^n$ (with the operator norm on $B^n$),
we have
\[
\|\Phi^n\| \le K
\]
for all $n$.
\end{lemma}

\begin{proof}
Since ${\mathfrak g}$ is finite dimensional, there is a constant, $J$,
such that $\|U_X\| \le J\|X\|$ for every $X \in {\mathfrak g}$.  Suppose
now that $X \in {\mathfrak g}_0$.  Then $U_X$ is skew-adjoint, and we can
find an orthonormal basis, $\{e_j\}$, for ${\mathcal H}$ consisting of
eigenvectors for $U_X$, with corresponding eigenvalues $\{\a_j\}$.  For
each $n$ we have the corresponding orthonormal basis $\{e_{j_1} \otimes
\dots \otimes e_{j_n}\}$ for ${\mathcal H}^{\otimes n}$.  Each $e_{j_1}
\otimes \dots \otimes e_{j_n}$ is an eigenvector for $U_X^{\otimes n}$,
with eigenvalue $\sum_{k=1}^n \a_{j_k}$.  (Note that here and below
$U_X^{\otimes n}$ denotes the {\em inner} tensor representation of
the Lie algebra ${\mathfrak g}$, 
not the tensor power of the operator $U_X$.)  Since
$\|U_X\|$ and $\|U_X^{\otimes n}\|$ are given by the largest absolute
value of the eigenvalues of $U_X$ and $U_X^{\otimes n}$, we see that
$\|U_X^{\otimes n}\| = n\|U_X\|$.  For $X$ in ${\mathfrak g}$ rather than
${\mathfrak g}_0$ it follows that $\|U_X^{\otimes n}\| \le 2n\|U_X\|$.
Upon restricting $U_X^{\otimes n}$ to ${\mathcal H}^n$ we then have
$\|U_X^n\| \le 2n\|U_X\| \le 2nJ\|X\|$.  Thus $\|\Phi_X^n\| =
\|n^{-1}U_X^n\| \le 2J\|X\|$, independent of $n$.

Now pick a basis for ${\mathfrak g}$.  For each $k$ it gives a basis for
${\mathfrak g}^{\otimes k}$ consisting of elementary tensors.  But for any
elementary tensor $Z = X_1 \otimes \dots \otimes X_k$ we have
\[
\|\Phi_Z^n\| = \|\Phi_{X_1}^n\dots \Phi_{X_k}^n\| \le (2J)^k \prod_{j=1}^k
\|X_j\|.
\]
The right-hand side is independent of $n$.  Since ${\mathcal T}^q$ is
finite dimensional the desired result follows easily.
\end{proof}

We let $\|\cdot\|_2$ denote the usual Hilbert-space norm on $L^2(G/H)$.
Let ${\mathcal T}^q_S$ denote the $S$-isotypic subspace of ${\mathcal
T}^q$.  Since $\Phi$ is equivariant, it carries ${\mathcal T}^q_S$ onto
$A_S$.

\setcounter{notation}{5}
\begin{notation}
\label{note6}
Let ${\mathcal F}$ $(= {\mathcal F}_S)$ denote the orthogonal complement
of the kernel of the restriction of $\Phi$ to ${\mathcal T}_S^q$.
\end{notation}

Thus $\Phi$ is a bijection from ${\mathcal F}$ onto $A_S$, and ${\mathcal
F}$ is carried into itself by $\Ad$.  Our notation will not distinguish
between $\Phi$ and its restriction to ${\mathcal F}$.  Much as above, we
let $B_S^n$ denote the $S$-isotypic subspace of $B^n$.

\setcounter{proposition}{6}
\begin{proposition}
\label{prop4.7}
There is an integer, $N$, such that for $n \ge N$ we have
\[
\s^n(\Phi^n({\mathcal F})) = A_S.
\]
In particular, $\s^n$ will for $n \ge N$ be a bijection from $B_S^n$ onto
$A_S$, and $\Phi^n$ will be a bijection from ${\mathcal F}$ onto $B_S^n$.
\end{proposition}

\begin{proof}
Put on $A_S$ the inner product from $L^2(G/H)$, and let $\{e_j\}$ be an
orthonormal basis for $A_S$.  For each $j$ let $Z_j$ be the unique element
of ${\mathcal F}$ such that $\Phi_{Z_j} = e_j$.  Let $d$ denote the
dimension of $A_S$.  All norms on a finite-dimensional vector space are
equivalent, and so from Theorem~\ref{th4.2} we see that we can find 
an integer $N$
such that
\[
\|\s^n(\Phi_{Z_j}^n) - \Phi_{Z_j}\|_2 < 1/d
\]
for all $n \ge N$ and all $j$.  Then for each $n \ge N$ the $f_j =
\s^n(\Phi_{Z_j}^n)$ span $A_S$.  To see this, write any $g \in A_S$ as $g
= \sum \a_je_j$, and set $h = \sum \a_jf_j$.  Then $\|g-h\|_2 \le \sum
|\a_j|/d \le \|g\|_2$, so that $g$ can not be orthogonal to the span of the
$f_j$'s unless it is 0.

Since $\s^n$ is injective, it must then be bijective 
from $B^n_S$ for $n \ge N$.  Since
the dimension of ${\mathcal F}$ is the same as that of $A_S$, it follows
that $\Phi^n$ is bijective from ${\mathcal F}$ onto $B^n_S$ for $n \ge
N$.
\end{proof}

We fix $N$ as above, and revert to the $C^*$-norms, restricted to $B_S^n$
and $A_S$.  Momentarily define $\Omega^n$ on $A_S$, for $n \ge N$, by
$\Omega^n = \s^n \circ \Phi^n \circ \Phi^{-1}$, where $\Phi^{-1}$ takes
$A_S$ to ${\mathcal F}$.  Note that each $\Omega^n$ is invertible.  Now
Theorem~\ref{th4.2} tells us that $\Omega^n(f)$ converges to $f$ for each
$f \in A_S$.  It follows that $\{\Omega^n\}$ converges in operator norm to
the identity operator.  Then so does $\{(\Omega^n)^{-1}\}$, and so this
sequence is uniformly bounded in norm.  Since $\Phi$ is independent of
$n$, we obtain:

\setcounter{lemma}{7}
\begin{lemma}
\label{lem4.8}
There is a constant, $r$, such that, as operators from $A_S$ to
${\mathcal F}$,
\[
\|(\s^n \circ \Phi^n)^{-1}\| \le r
\]
for all $n \ge N$.
\end{lemma}

\setcounter{corollary}{8}
\begin{corollary}
\label{cor4.9}
There is a constant, $K'$, such that for each $n \ge N$, when we view
$(\s^n)^{-1}$ as an operator from $A_S$ to $B_S$, we have
\[
\|(\s^n)^{-1}\| \le K'.
\]
\end{corollary}

\begin{proof}
 From the above lemma and from Lemma~\ref{lem4.5} we have
\[
\|(\s^n)^{-1}\| = \|\Phi^n \circ (\s^n \circ \Phi^n)^{-1}\| \le Kr.
\]
\end{proof}

For $r$ as above 
let ${\mathcal B}_r$ denote the closed ball of radius $r$ in ${\mathcal
F}$.  It follows from Lemma~\ref{lem4.8} that every $f$ in the unit ball
of $A_S$ is, for every $n \ge N$, of the form $\s^n(\Phi_{Z_n}^n)$ for
some $Z_n \in {\mathcal B}_r$.  Let $T \in B_S^n$ with $\|T\| \le 1$.
Then $\|\s^n(T)\| \le 1$, and so if $n \ge N$ we conclude from the above
observation that $T = \Phi_Z^n$ for some $Z \in {\mathcal B}_r$.  We have
thus obtained (for $N$ as in Proposition~\ref{prop4.7}):

\setcounter{proposition}{9}
\begin{proposition}
\label{prop4.10}
There is a closed ball, ${\mathcal B}$, in ${\mathcal F}$ such that
$$
\Phi^n({\mathcal B}) \supseteq (\mbox{\em unit ball of } B_S^n)
$$
for every $n \ge N$.
\end{proposition}


\section{Contravariant symbols}
\label{sec5}

This section is devoted to obtaining the information about ${\breve \s}^n$
which we need.  The essence of what we need is contained in theorem 2B of
\cite{Dff}.  This attractive theorem, concerning convergence of
contravariant symbols as $n$ increases, seems to have few antecedents in
the literature, and seems not to have been used at all since its
appearance.  But it is crucial to our present purposes.  However, 
we need tighter control over the situation than is
explicitly given in Duffield's proof.  Part of this section will basically
consist of rewriting Duffield's proof so as to give this tighter control
(and with simpler arguments).  We will also draw some important
consequences.

We continue with the notation of the previous section.  Because ${\breve
\s}^n$ is equivariant, it will carry $A_S$ into $B_S^n$.  We continue to
let the integer $N$ be as in Proposition~\ref{prop4.7}, so that $\s^n$ is
a bijection from $B_S^n$ onto $A_S$ for $n \ge N$.  Then ${\breve \s}^n$
will be a bijection from $A_S$ onto $B_S^n$ for $n \ge N$, and $({\breve
\s}^n)^{-1}$ will exist from $B_S^n$ to $A_S$.  When we write $({\breve
\s}^n)^{-1}$ in the next pages, we will always take it as defined on
$B_S^n$ for $n \ge N$.  Note that if $T \in B_S^n$ and if we set $f =
({\breve \s}^n)^{-1}(T)$, then ${\breve \s}_f^n = T$.  Thus $({\breve
\s}^n)^{-1}$ provides a canonical way of choosing a contravariant symbol
for $T$.  Recall that 
for $n \ge N$ each $\Phi^n$ is a bijection from ${\mathcal F}$
onto $B_S^n$ according to Proposition~\ref{prop4.7}.

\begin{notation}
\label{note5.1}
For $n \ge N$ set $\Psi^n = ({\breve \s}^n)^{-1} \circ \Phi^n$ on ${\mathcal
F}$, so that $\Psi^n$ is a bijection from ${\mathcal F}$ onto $A_S$.
\end{notation}

Thus $\Psi_Z^n$ is a contravariant symbol for $\Phi_Z^n$ for each $Z \in
{\mathcal F}$.  We need to show, in essence, that for each $Z \in
{\mathcal F}$ the $\Psi_Z^n$'s converge to $\Phi_Z$.

For $n \ge N$ define a linear functional, $\th_n$, on ${\mathcal F}$ by
\[
\th_n(Z) = \Psi_Z^n(e).
\]
Because $\Psi_Z^n \in C(G/H)$, we have $\th_n(\Ad_s(Z)) = \th_n(Z)$ for
every $s \in H$.  Because $\Psi^n$ is equivariant, we have
\[
\th_n(\Ad_x(Z)) = \Psi^n_{\Ad_x(Z)}(e) = 
(\l_x(\Psi^n_Z))(e) = \Psi_Z^n(x^{-1})
\]
for every $x \in G$.  Since $\Psi_Z^n$ is a contravariant symbol for
$\Phi_Z^n$, we obtain, as seen early in Section~\ref{sec2},
\[
\Phi_Z^n = d_n \int_{G/H} \th_n(\Ad_{y^{-1}}(Z))\a_y(P^n)dy.
\]
 From formula \eqref{eq2.2}, together with the notation $h^n = h_{P^n}$
used in Section~\ref{sec3}, we then obtain
\[
\s^n(\Phi_Z^n)(x) = \int_{G/H} \th_n(\Ad_{yx^{-1}}(Z))h^n(y)dy.
\]
Now Theorem~\ref{th4.2} tells us that the left-hand side converges to
$\Phi_Z(x)$.  In view of Lemma~\ref{lem3.3}, this suggests that perhaps
$\th_n(\Ad_{x^{-1}}(Z))$ also converges to $\Phi_Z(x)$.  But to show that
this is correct we need control of the size of the $\th_n$'s.  The key
fact which we will use for this is Lemma~\ref{lem4.5}, which tells us that
$\|\Phi^n\| \le K$ for all $n$.

We now use the inner product  on ${\mathcal F} \subseteq {\mathcal T}^q$.
Then for each $n \ge N$ there will be a $Z_n \in {\mathcal F}$ which
represents $\th_n$, so that $\th_n(Z) = \<Z,Z_n\>$ for all $Z \in
{\mathcal F}$.  Then $\|\th_n\| = \|Z_n\|$, where we use the norm from the
inner product.  Because $\th_n(\Ad_s(Z)) = \th_n(Z)$ for every $s \in H$,
we see quickly that $Z_n$ is $H$-invariant (for $\Ad$).

Let ${\mathcal F}^H$ denote the subspace of $H$-invariant elements of
${\mathcal F}$, so that each $Z_n \in {\mathcal F}^H$.  As $Z$ ranges
over ${\mathcal F}^H$ the functions $y \mapsto \Ad_y(Z)$ form a
finite-dimensional vector space of vector-valued functions on $G/H$.  By
considering a basis for this vector space, we see that any bounded
collection of these functions will be equicontinuous.  From
Lemma~\ref{lem3.3} we find that for every $\e > 0$ there is an integer
$M_{\e} \ge N$ such that
\[
\left\| Z - \int_{G/H} \Ad_y(Z)h^n(y)dy\right\| \le \e\|Z\|
\]
for all $Z \in {\mathcal F}^H$ and all $n \ge M_{\e}$.  Then for $n \ge
M_{\e}$ we have
\begin{eqnarray*}
\|Z_n||^2 &= &\<Z_n,Z_n - \int \Ad_y(Z_n)h^n(y)dy\> + \<Z_n, \int
\Ad_y(Z_n)h^n(y)dy\> \\
&\le &\e\|Z_n\|^2 + \left\vert \int \th_n(\Ad_y(Z_n))h^n(y)dy\right\vert.
\end{eqnarray*}
Upon applying Lemma~\ref{lem4.5}, with $K$ as given there, we then obtain
\[
(1-\e)\|Z_n\|^2 \le |\s^n(\Phi_{Z_n}^n)(e)| \le \|\Phi^n\|\|Z_n\| \le
K\|Z_n\|.
\]
Since $\|Z_n\|=\|\th_n\|$, this gives the proof of:

\setcounter{lemma}{1}
\begin{lemma}
\label{lem5.2}
For every $r > 1$ there is an integer $M_r \geq N$ such that for $n \ge M_r$ we
have $\|\th_n\| \le rK$.
\end{lemma}

 From this we see that for $n \geq N$ and any $Z \in \mathcal F$ we have
\[
|\Psi_Z^n(x)| = |\th_n(\Ad_{x^{-1}}(Z))| \le rK\|Z\|,
\]
for all $x$, from which we obtain:

\setcounter{proposition}{2}
\begin{proposition}
\label{prop5.3}
For every $r > 1$ there is an integer $M_r \geq N$ such that 
if $n \ge M_r$ then
$\|\Psi^n\| \le rK$ (where $\Psi^n$ is defined on ${\mathcal F}$, and both
${\mathcal F}$ and $K$ depend on $S$).
\end{proposition}

The following theorem is our version of theorem 2B of \cite{Dff}.

\setcounter{theorem}{3}
\begin{theorem}
\label{th5.4}
For any $\e > 0$ there is an integer $N_{\e} \geq N$ 
such that if $n \ge N_{\e}$
then
\[
\|\Psi_Z^n - \s^n(\Phi_Z^n)\|_{\infty} \le \e\|Z\|
\]
for all $Z \in {\mathcal F}$ ($={\mathcal F}_S)$.  
In particular, $\Psi_Z^n$ converges to
$\Phi_Z$.
\end{theorem}

\begin{proof}
For any $Z \in {\mathcal F}$ we have, much as above,
\[
\Psi_Z^n(e) - \s^n(\Phi_Z^n)(e) = \int_{G/H} \th_n(Z -
\Ad_{y^{-1}}(Z))h^n(y)dy.
\]
As $Z$ ranges over ${\mathcal F}$ the functions $y \mapsto Z -
\Ad_{y^{-1}}(Z)$ form a finite-dimensional vector space of vector-valued
functions on $G$ (not $G/H$), so, as before, any bounded subset will be
equicontinuous.  We know from Lemma~\ref{lem5.2} that $\{\|\th_n\|: n \ge
N\}$ is bounded.  Thus the set of functions $y \mapsto \th_n(Z -
\Ad_{y^{-1}}(Z))$, now on $G/H$, is equicontinuous for $\|Z\| \le 1$ and
for all $n \ge N$.  Let $\e > 0$ be given.  It follows from
Lemma~\ref{lem3.3} that we can find an integer, $N_{\e}$, such that
\[
|\Psi_Z^n(e) - \s^n(\Phi_Z^n)(e)| \le \e\|Z\|
\]
for all $n \ge N_{\e}$ and all $Z \in {\mathcal F}$.  But from the
equivariance of $\Psi^n$ and $\s^n \circ \Phi^n$, and the fact that $\Ad$
is unitary on ${\mathcal F}$, it follows that
\[
\|\Phi_Z^n - \s^n(\Phi_Z^n)\|_{\infty} \le \e\|Z\|
\]
for $n \ge N_{\e}$ and $Z \in {\mathcal F}$.
\end{proof}


\section{Conclusion of the proof of Theorem \ref{th3.2}}
\label{sec6}

We continue with the notation of the previous sections.  Our first
objective is to prove the following fact, which seems to be new and of
independent interest.

\begin{theorem}
\label{th6.1}
Let the general hypotheses be the same as those for Theorem 3.4.
For each $n \ge 1$ let $\g_n$ be the smallest constant such that
\[
\|T - {\breve \s}^n(\s_T^n)\| \le \g_nL_n(T)
\]
for all $T \in B^n$.  Then the sequence $\{\g_n\}$ converges to $0$.
\end{theorem}

We remark that our wording here is somewhat different from
that of Theorem 3.4 because in the present finite-dimensional
situation it is clear that the constants $\g_n$ exist, whereas
in Theorem 3.4 it is not immediately clear that (finite)
constants $\g_n$ exist.

\begin{proof}[Proof of Theorem \ref{th6.1}]
Let $\e > 0$ be given.  By theorem $8.2$ and lemma $8.3$ of \cite{R6} we
can find a finite subset $S \subseteq {\hat G}$ and a positive linear 
combination, $\varphi$, of the
characters of the elements of $S$, with the following properties.
The set $S$ contains the
trivial representation $1$ and is closed under taking contragradient
representations; and for any ergodic action $\a$ of
$G$ on a unital $C^*$-algebra $C$ the integrated operator $\a_{\varphi}$
is a completely positive unital equivariant map of $C$ onto its
$S$-isotypic component such that $\|c - \a_{\varphi}(c)\| \le (\e/3)L(c)$
for all $c \in C$.  Then for every $T \in B^n$ we have
\[
\|T - {\breve \s}^n(\s_T^n)\| \le (\e/3)L_n(T) + \|\a_{\varphi}(T) -
{\breve \s}^n(\s_{\a_{\varphi}(T)}^n)\| + (\e/3)L_n(T).
\]
Thus we see that it suffices to prove that, for $S$ fixed as above, there is
an integer $N$ such that for all $T \in B_S^n$ we have
\[
\|T - {\breve \s}^n(\s_T^n)\| \le (\e/3)L_n(T)
\]
for all $n \ge N$.

Now by lemma $2.4$ of \cite{R4} all of the algebras $B^n$ will have radius
no larger than $r = \int_G \ell(x)dx$, in the sense that $\|T\|^{\sim} \le
rL(T)$ for all $T$, where $\|\cdot\|^{\sim}$ denotes the quotient of
$\|\cdot\|$ on $B^n/{\mathbb C}I$.  It then suffices to prove:

\setcounter{lemma}{1}
\begin{lemma}
\label{lem6.2}
For $S$ fixed as above, there is an integer $N$ such that
\[
\|T - {\breve \s}^n(\s_T^n)\| \le (\e/3r)\|T\|
\]
for all $n \ge N$ and all $T \in B_S^n$.
\end{lemma}

To see that this is sufficient, note that if it holds, then it holds
equally well for $T + tI$ where $t \in {\mathbb C}$ is chosen so that
$\|T\|^{\sim} = \|T + tI\|$, and that the left side of the inequality will
be unchanged.

\begin{proof}[Proof of Lemma \ref{lem6.2}]
We use Lemma~\ref{lem4.4} and Notation 4.6 
to choose $q$ and ${\mathcal F} \subset
{\mathcal T}^q$ for our fixed $S$.  We then apply Proposition~\ref{prop4.10} to
choose an integer $N$ such that there is a closed ball ${\mathcal B}$ in
${\mathcal F}$ such that $\Phi^n({\mathcal B}) \supseteq (\mbox{unit ball
of } B_S^n)$ for all $n \ge N$.  Let $R$ denote the radius of ${\mathcal
B}$, and let $\e' = \e/3rR$.  We apply Theorem~\ref{th5.4} to choose a yet
larger  $N$ such that $\Psi^n_Z$ is defined for every $n \ge N$, 
and for all $Z \in {\mathcal F}$ we have
\[
\|\Psi_Z^n - \s^n(\Phi_Z^n)\|_{\infty} \le \e'\|Z\|   .
\]

Suppose now that $n \ge N$ and that $T \in B_S^n$ with $\|T\| \le 1$.  By
our choice of $N$ there is a $Z \in {\mathcal B}$ such that $\Phi_Z^n =
T$, and furthermore $T = {\breve \s}^n(\Psi_Z^n)$.  Then
\begin{eqnarray*}
\|T - {\breve \s}^n(\s_T^n)\| &= &\|{\breve \s}^n(\Psi_Z^n) - {\breve
\s}^n(\s_T^n)\| \\
&\le &\|\Psi_Z^n - \s^n(\Phi_Z^n)\|_{\infty} \le \e'\|Z\| \le \e/3r.
\end{eqnarray*}
Since this is true for all $T$ with $\|T\| \le 1$, we obtain the desired
result.
\end{proof}

We have thus concluded the proof of Theorem~\ref{th6.1}.

We now complete the proof of our main theorem, Theorem~\ref{th3.2}.  Let
$\e > 0$ be given.  By Theorem~\ref{th3.4} we can find an integer $N$ such
that for $n \geq N$ we have
\[
\|f - \s^n({\breve \s}_f^n)\|_{\infty} \le (\e/2)L_A(f)
\]
for all $f \in A$.  By Proposition~\ref{prop1.3} and
Corollary~\ref{cor2.4}, it follows that $S(A)$ is in the 
$(\e/2)$-neighborhood of $S(B^n)$ for $\rho_L$, for each $n \ge N$.

Now according to Theorem~\ref{th6.1} we can find $N$ still larger such
that for $n \ge N$ we have
\[
\|T - {\breve \s}^n(\s_T^n)\| \le (\e/2)L_n(T)
\]
for $T \in B^n$.  Then by the calculation done in the next-to-last
paragraph of Section~\ref{sec2} (where the $\g$ there is the $\e/2$ here), 
it follows that each $S(B^n)$ is
in the $\e$-neighborhood of $S(A)$ for $\rho_L$.  Accordingly, for $n \ge N$
\[
\dist_q((B^n,L_n),(A,L_A)) \le \e,
\]
as desired.
\end{proof}

We remark that when one examines the steps in the proof of our main
theorem, one sees that for any given coadjoint orbit and any given length
function $\ell$ one can with careful bookkeeping obtain, for any $\e > 0$,
a fairly explicit choice of an $N$ for which
\[
\dist_q((B^n,L_n),(A,L_A)) \le \e
\]
for $n \ge N$.


\end{document}